\documentclass[12pt,draftcls,onecolumn]{IEEEtran}

\usepackage{times} 
\usepackage{amsmath}
\usepackage{amsfonts}
\usepackage{amssymb}
\usepackage{graphicx}
\DeclareGraphicsExtensions{.pdf,.jpeg,.png}
\usepackage{pdfpages}
\usepackage{epstopdf}

\title{\bf Observability Analysis of Sensorless Synchronous Machine Drives}

\author{Mohamad Koteich, \textit{Student Member, IEEE}, \\ Abdelmalek Maloum, Gilles Duc and Guillaume Sandou
	\thanks{Mohamad Koteich is with Renault S.A.S. Technocentre, 78288 Guyancourt, France, and also with L2S - CentraleSup\'{e}lec - CNRS - Paris-Sud University, 91192 Gif-sur-Yvette, France (e-mail: mohamad.koteich@renault.com).}
	\thanks{Abdelmalek Maloum is with Renault S.A.S. Technocentre, 78288 Guyancourt, France (e-mail: abdelmalek.maloum@renault.com).}
	\thanks{Gilles Duc and Guillaume Sandou are with L2S - CentraleSup\'{e}lec - CNRS - Paris-Sud University, 91192 Gif-sur-Yvette, France  (e-mail: gilles.duc@centralesupelec.fr; guillaume.sandou@centralesupelec.fr).}}

\begin{document}
\maketitle

\begin{abstract}
This paper studies the local observability of synchronous machines using a unified approach. Recently, motion sensorless control of electrical drives has gained high interest. The main challenge for such a technology is the poor performance in some operation conditions. One interesting theory that helps understanding the origin of this problem is the observability analysis of nonlinear systems. In this paper, the observability of the wound-rotor synchronous machine is studied. The results are extended to other synchronous machines, adopting a unified analysis. Furthermore, a high-frequency injection-based technique is proposed to enhance the sensorless operation of the wound-rotor synchronous machine at standstill.

\end{abstract}
\section{Introduction}

Electrical rotating machines are becoming very popular in nowadays transport industries, such as electric vehicle and more electric aircraft \cite{emadi14} \cite{mea12}. 

Synchronous machine (SM) is one of the biggest families of electrical machines, which is widely used in high performance industry applications. Various types of SMs can be classified depending on the rotor configuration \cite{chiasson}; there exist wound-rotor (WRSM), permanent-magnet (PMSM) and reluctance type (SyRM) synchronous machines.

Over the last few decades, many control techniques have been proposed and used for electrical drives \cite{boldea08}. Mechanical sensorless techniques \cite{vas} \cite{holtz06} \cite{koteich13} have been good candidates for reliable and costless electrical drives \cite{betin14}. Nevertheless, these techniques have the problem of deteriorated performance in some operation conditions. 

Recently, observability analysis of electrical drives, based on local weak observability theory of nonlinear systems \cite{hermann77}, has taken more interest in order to understand observer's deteriorated performance. 

In contrast to observability of linear systems, observability of nonlinear systems depends on the inputs and initial conditions. An observable nonlinear system might be unobservable with some inputs (singular inputs), which affects the observer operation \cite{besancon}.

Observability of induction machines (IM) is studied in \cite{canudas00}, \cite{rojas} and \cite{chiasson06}. More recently, the observability study of SMs has started only for the PMSM \cite{zaltni10} \cite{ezzat10} \cite{abry11}. To the best of the authors knowledge, the first paper that could formulate useful explicit observability conditions for the PMSM is \cite{vaclavek13}, where the conditions are expressed in the rotor reference frame.

In the present work, the WRSM observability is analyzed, and the results are extended to the other SMs using a unified approach. Furthermore, based on the aforementioned analysis, a high-frequency (HF) injection-based technique is proposed, in order to ensure the WRSM observability in the unobservable region. The results are validated using an Extended Kalman Filter (EKF) and illustrated via numerical simulations.

The main result of the unified observability analysis is the definition of a fictitious \textit{observability vector} for SMs; the local observability of any SM is guaranteed as soon as the rotational velocity of the observability vector with respect to the rotor is different from the electrical velocity of the rotor with respect to the stator.

This paper is organized as follows: in section II, the local observability concept of nonlinear systems is presented. In section III, the state-space model of the WRSM is given, the other SMs models are derived from the WRSM one. The observability of SMs is studied in section IV. Section V presents illustrative simulations that validate the obtained results with the proposed HF injection technique.

\section{Observability theory}
\label{obsv_theory}
There are many approaches to study the observability of nonlinear systems. In this section, the \emph{local weak observability} concept \cite{hermann77}, based on the rank criterion approach, is presented. This approach provides only sufficient conditions.

\subsection{Problem statement}
Systems of the following form (denoted $\Sigma$) are considered:
\begin{equation}
\Sigma :\left\{
\begin{aligned}
\dot{x} &= f\left(x(t), u(t)\right)\\
y &= h\left(x(t)\right)
\label{sigma}
\end{aligned}
\right.\end{equation}
where $x \in X \subset \mathbb{R}^n$ is the state vector, $u \in U \subset \mathbb{R}^m$ is the control vector (input), $y \in \mathbb{R}^p $ is the output vector, $f$ and $h$ are $C^\infty$ functions.

The observation problem can be then formulated as follows \cite{besancon}:
\emph{Given a system described by a representation \eqref{sigma}, find an accurate estimate $\hat{x}(t)$ for $x(t)$ from the knowledge of $u(\tau)$, $y(\tau)$ for $0 \leq \tau \leq t$}.

\subsection{Definitions}
\paragraph{Indistinguishability}
Let $x_0$ and $x_1$ be two initial states of the system $\Sigma$ \eqref{sigma} at the time $t_0$ ($x_0, x_1 \in x \subset X$). The pair $(x_0, x_1)$ is indistinguishable if, for any admissible input $u(t)$, the system outputs $y_0(t)$ and $y_1(t)$, respectively associated to $x_0$ and $x_1$, follow the same trajectories from $t_0$ to $t$, i.e. starting from those two initial states, the system realizes the same input-output map \cite{hermann77}. Otherwise, $x_0$ and $x_1$ are distinguishable.

\paragraph{Observability}
A system \eqref{sigma} is observable (resp. at $x_0$) if it does not admit any indistinguishable pair (resp. any state indistinguishable from $x_0$) \cite{besancon}.

This definition is too general. In practice, one might be interested in distinguishing states from their neighbors.

\paragraph{Local weak observability}
A system \eqref{sigma} is locally weakly observable (resp. at $x_0$) if there exists a neighborhood $V$ of any $x$ (resp. of $x_0$) such that for any neighborhood $W$ of $x$ (resp. $x_0$) contained in $V$, there is no indistinguishable state from $x$ (resp. $x_0$) in $W$ when considering time intervals for which trajectories remain in $W$.
This roughly means that one can distinguish every state from its neighbors without ``going too far''. This notion is of more interest in practice, and also presents the advantage of admitting some `rank condition' characterization \cite{besancon}.

\paragraph{Observation space} The observation space for a system \eqref{sigma} is defined as the smallest real vector space (denoted by $\mathcal{S}_h$) of $C^\infty$ functions containing the components of $h$ and closed under \emph{Lie derivation} along $f_u:= f(., u)$ for any constant $u \in \mathbb{R}^m$:
\begin{equation}
\mathcal{S}_h(x) = \left[ \begin{matrix}
\mathcal{L}^0_fh(x)\\
\mathcal{L}_fh(x)\\
\mathcal{L}_f^2h(x)\\
\vdots\\
\mathcal{L}_f^{n-1}h(x)
\end{matrix} \right]
\end{equation}
where $\mathcal{L}_f^{k}h$ is the $k$th-order \emph{Lie derivative} of the function $h$ with respect to the vector field $f$. 

\subsection{Observability rank condition}
The system $\Sigma$ is said to satisfy the observability rank condition at $x_0$ if the Jacobian of the observability space (called observability matrix and denoted by $\mathcal{O}_h(x)$) is full rank at $x_0$:
\begin{equation}
 rank~\mathcal{O}_h(x)|_{x_0} = rank \left[\frac{\partial\mathcal{S}_h(x)}{\partial x}\right]_{x_0} = n
\end{equation}

\subsection{Observability theorem}
From the previous definitions, the following theorem can be stated \cite{hermann77}:
A system $\Sigma$ \eqref{sigma} satisfying the observability rank condition at $x_0$ is locally weakly observable at $x_0$. More generally, a system $\Sigma$ \eqref{sigma} satisfying the observability rank condition, for any $x_0$, is locally weakly observable.

\section{Synchronous machine model}
\label{ss_model}

In this section, the mathematical model of the WRSM is presented. The models of other SMs can be extended from the WRSM one. The assumption of linear lossless magnetic circuit is adopted, with sinusoidal distribution of stator magnetomotive force. The machine parameters are considered to be known constants. Nevertheless, the parameters variation does not call the observability study results into question; it impacts the observer performance, which is beyond the scope of this study.

\subsection{Machine description}
\begin{figure}[!t]
\centering
\includegraphics[scale=0.45]{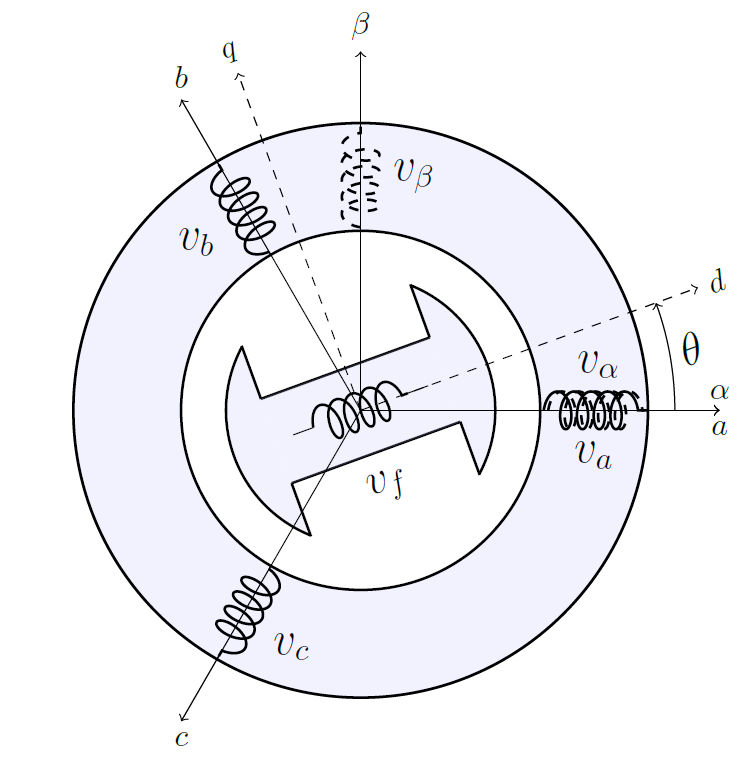}
\caption{Symbolic representation of the salient-type wound-rotor synchronous machine}
\label{wrsm_ab}
\end{figure}

Synchronous machines are electromechanical systems composed of two parts (see Fig. \ref{wrsm_ab}):
\begin{itemize}
\item Stator, the stationary part, fed by a three-phase source.
\item Rotor, the moving part, which defines the sub-family of an SM depending on its type:
\begin{itemize}
\item[1-] WRSM: the rotor is an electromagnet supplied by a DC source.
\item[2-] PMSM: the rotor is made of permanent magnets that can be Interior (IPMSM), or Surface-mounted (SPMSM).
\item[3-] SyRM: the rotor has neither permanent magnets nor windings, it is made of a ferromagnetic core.
\end{itemize}
\end{itemize}
Both the WRSM and PMSM can be either salient-type (non cylindrical) rotor, that is airgap between stator and rotor varies as the rotor moves, or non-salient type (cylindrical) rotor. As for the SyRM, its rotor is necessarily salient type, since the operation principle of this machine is based on rotor alignment with the stator rotating magnetic field following the minimum reluctance magnetic path.

\subsection{State-space model of salient-type WRSM}
The electromagnetic behaviour of the stator windings in a three-phase electrical machine, such as the WRSM, can be fully described using two equivalent (fictitious) two-phase stator windings \cite{park29} \cite{park33}, denoted $\alpha$ and $\beta$ (Fig. \ref{wrsm_ab}). The state-space model of the salient-type WRSM can be written, in the ($\alpha \beta$) stationary reference frame, in a way to be fitted to the structure \eqref{sigma}:
\begin{eqnarray}
\label{ss_ab_wrsm}
\dot{x}&=& f(x, u) : \left\{
\begin{aligned}
\frac{d{\mathcal{I}}}{dt}&=-{\mathfrak{L}^{-1}}{\mathfrak{R}_{eq}}\mathcal{I}
+ {\mathfrak{L}^{-1}}\mathcal{V}\\
\frac{d\omega}{dt} &= - \frac{f_v}{J}\omega + \frac{p}{J} T_m - \frac{p}{J} T_l\\
\frac{d\theta}{dt} &= \omega
\end{aligned}\right. 
\label{ss}\\
y &=& h(x) = \mathcal{I}
\label{out}
\end{eqnarray}
where the state, input and output vectors are respectively:
\begin{eqnarray}
x = \left[
\begin{matrix}
\mathcal{I}^T & \omega & \theta
\end{matrix}
\right]^T;~
u = \mathcal{V}~;~
y = \mathcal{I}
\end{eqnarray}
The first equation in the system \eqref{ss} comes from the Ohm's law, where $\mathcal{I}$ and $\mathcal{V}$ are the current and voltage vectors:
\begin{eqnarray}
\mathcal{I} = \left[\begin{matrix}
i_\alpha & 
i_\beta & 
i_f
\end{matrix}\right]^T~~;~~
\mathcal{V} = \left[\begin{matrix}
v_\alpha & v_\beta & v_f
\end{matrix}\right]^T
\end{eqnarray}
Indices $\alpha$ and $\beta$ stand for stator signals, index $f$ stands for rotor (\textit{field}) ones. 

$\mathfrak{L}$ is the ($\theta$-dependent) matrix of inductances:
\begin{equation}
\mathfrak{L}=\left[ \begin{matrix}
   {{L}_{0}}+{{L}_{2}}\cos 2\theta  & {{L}_{2}}\sin 2\theta  & {{M}_{f}}\cos \theta   \\
   {{L}_{2}}\sin 2\theta  & {{L}_{0}}-{{L}_{2}}\cos 2\theta  & {{M}_{f}}\sin \theta   \\
   {{M}_{f}}\cos \theta  & {{M}_{f}}\sin \theta  & {{L}_{f}}  \\
\end{matrix} \right]
\end{equation}
where $L_0 = (L_d + L_q)/2$ and $L_2 = (L_d - L_q)/2$. $L_d$ and $L_q$ being the direct and quadrature inductances of the equivalent machine model in the rotor ($dq$) reference frame (Fig. \ref{wrsm_ab}) \cite{park29} \cite{park33}. $L_f$ is the rotor winding inductance and $M_f$ is the maximal mutual inductance between stator and rotor windings.

$\mathfrak{R}_{eq}$ is the \textit{equivalent resistance} matrix defined as:
\begin{equation}
\mathfrak{R}_{eq} = \mathfrak{R} + \frac{\partial\mathfrak{L}}{\partial\theta}\omega
\end{equation}
$\mathfrak{R}$ is the matrix of resistances ($R_s$ and $R_f$ stand respectively for stator and rotor resistance):
\begin{equation}
\mathfrak{R} = \left[\begin{matrix}
R_s & 0   &  0 \\
  0 & R_s &  0 \\
  0 & 0   & R_f
\end{matrix}\right]
\end{equation}

$\omega$ denotes the electrical speed (rad/sec) and $\theta$ the electrical position of the rotor\footnote{electrical speed (resp. position) = $p~\times$ mechanical speed (resp. position)}.

The second state equation of the system \eqref{ss} comes from the Newton's second law for rotational motion, where $J$ is the moment of inertia of the rotor with its associated load, $f_v$ is the viscous friction coefficient, $p$ is the number of pole pairs, $T_l$ is the load torque and $T_m$ is the motor torque given by:
\begin{eqnarray}
T_m = \frac{3}{2}.\frac{p}{2}~ \mathcal{I}^T~ \frac{\partial \mathfrak{L}}{\partial \theta} ~\mathcal{I} 
\end{eqnarray}

\subsection{State-space model of the other SMs}
\label{subsec:ss_sms}
\begin{figure}
\centering
\includegraphics[scale=0.45]{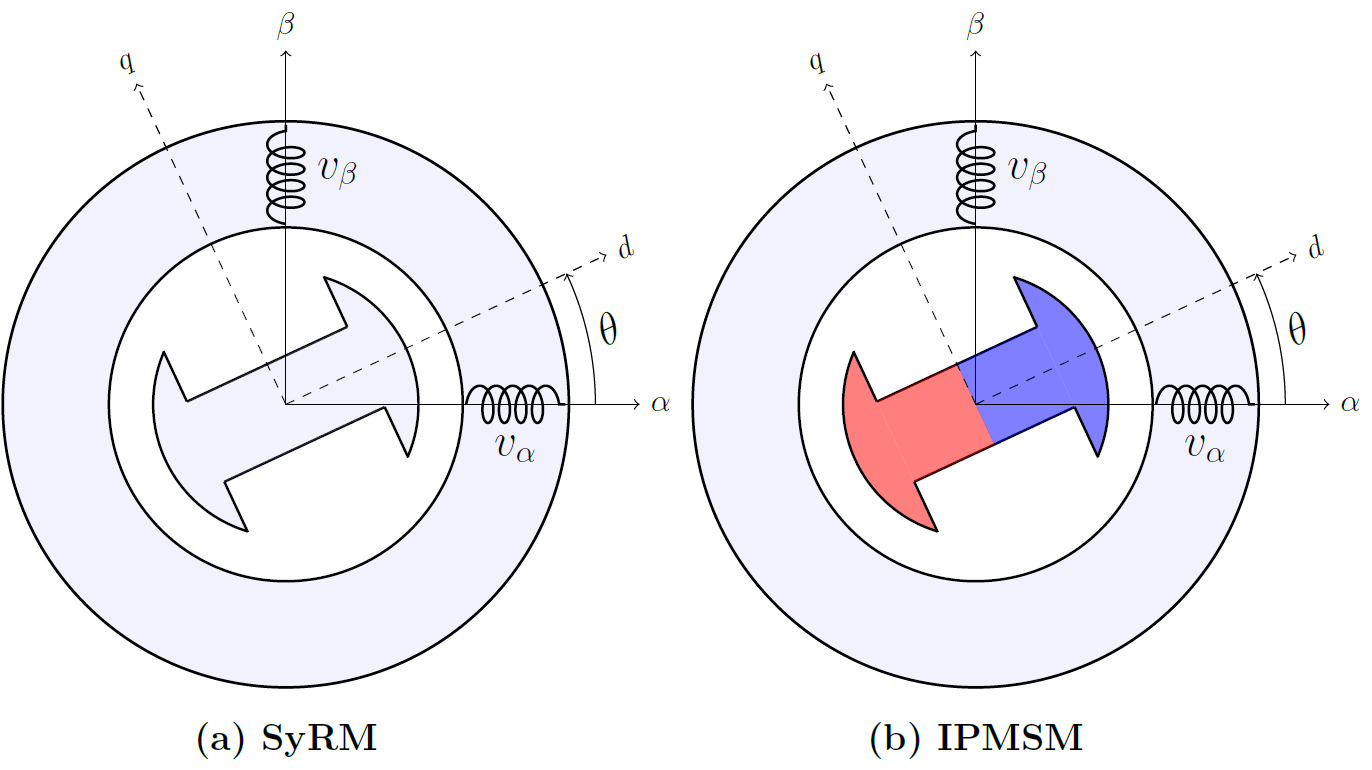}
\caption{Symbolic representation of SyRM(a) and IPMSM(b) in the $\alpha \beta$ reference frame}
\label{sm_ab}
\end{figure}

\begin{figure}
\centering
\includegraphics[scale=0.45]{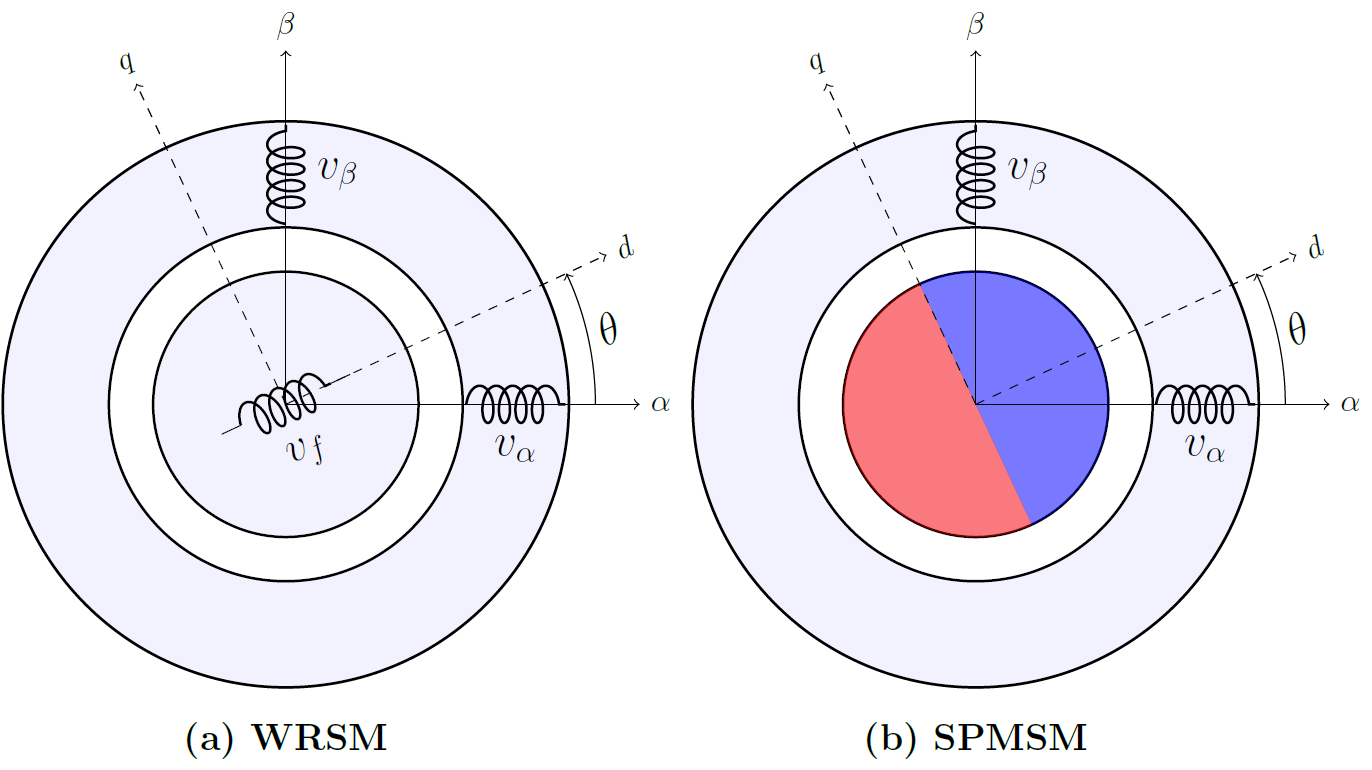}
\caption{Symbolic representation of non-salient WRSM(a) and SPMSM(b) in the $\alpha \beta$ reference frame}
\label{nsp_sm_ab}
\end{figure}

The other SMs can be seen as special cases of the salient-type WRSM; the IPMSM model (Fig. \ref{sm_ab}(b)) can be derived by considering the rotor magnetic flux to be constant:
\begin{eqnarray}
\frac{di_f}{dt} = 0
\label{cond_pmsm_1}
\end{eqnarray}
and by substituting $M_f i_f$ by the permanent magnet flux $\psi_r$:
\begin{eqnarray}
i_f = \frac{\psi_r}{M_f}
\label{cond_pmsm_2}
\end{eqnarray}
The SyRM model (Fig. \ref{sm_ab}(a)) can be derived from the IPMSM model by considering the rotor magnetic flux $\psi_r$ to be zero: 
\begin{eqnarray}
\psi_r \equiv 0
\label{cond_syrm}
\end{eqnarray}
The equations of the non-salient WRSM and SPMSM (Fig.~\ref{nsp_sm_ab}) are the same as the salient WRSM and IPMSM respectively, except that the stator self-inductances are constant and independent of the rotor position, that is:
\begin{eqnarray}
L_2 = 0 ~~~\implies~~~ L_d = L_q = L_0
\label{non_salient}
\end{eqnarray}

\section{Synchronous machine observability}
\label{obsv_analysis}
In this section, the local observability of the system \eqref{ss_ab_wrsm} is analyzed, in order to know if the mechanical states $\omega$ and $\theta$  can be estimated when only the currents $\mathcal{I}$ and voltages $\mathcal{V}$ are known. The state-space model \eqref{ss} is considered and the observability theory presented in section \ref{obsv_theory} is used. The machine model is strongly nonlinear; some calculations will be done using symbolic math software. It should be kept in mind that the observability rank condition is only a sufficient condition. 

\subsection{Observability matrix}
The system \eqref{ss} is a $5$-th order system. Its observability matrix should contain the output and its derivatives up to the $4$-th order. In this study, only the first order derivatives are calculated, higher order derivatives are very difficult to calculate and to deal with. This gives the following ``partial'' observability matrix:
\begin{equation}
\mathcal{O}_{y} =\left[ \begin{matrix}
\mathbb{I}_{3 \times 3} &\mathbb{O}_{3 \times 1} &\mathbb{O}_{3 \times 1}  \\
-\mathfrak{L}^{-1}\mathfrak{R}_{eq} &-\mathfrak{L}^{-1}\mathfrak{L}' \mathcal{I} &\mathfrak{L^{-1}}' \mathfrak{L} \frac{d \mathcal{I}}{dt} - \mathfrak{L^{-1}}\mathfrak{L}''\omega \mathcal{I} 
\end{matrix} \right]
\label{obsv_matrix}
\end{equation}
where $\mathbb{I}_{n \times n}$ is an $n \times n$  identity matrix, and $\mathbb{O}_{n \times m}$ is an $n \times m$ zero matrix. $\mathfrak{L}'$ and $\mathfrak{L}''$ denote, respectively, the first and second partial derivatives of $\mathfrak{L}$ with respect to $\theta$:
\begin{eqnarray}
\mathfrak{L}' = \frac{\partial}{\partial \theta} \mathfrak{L}~~~;~~
\mathfrak{L}'' = \frac{\partial}{\partial \theta} \mathfrak{L}'
\end{eqnarray}
The matrix \eqref{obsv_matrix} is a $6\times5$ matrix. It is full-rank if, at least, one of its $5\times5$ sub-matrices is full-rank. Regarding the structure of the matrix \eqref{obsv_matrix}, the rank study can be restricted to the following $3\times2$ sub-matrix:
\begin{equation}
\begin{bmatrix}
-\mathfrak{L}^{-1}\mathfrak{L}' \mathcal{I} &,&  \mathfrak{L^{-1}}' \mathfrak{L} \frac{d \mathcal{I}}{dt} - \mathfrak{L^{-1}}\mathfrak{L}''\omega \mathcal{I} 
\end{bmatrix}
\label{sub_matrix_obsv}
\end{equation}
It is sufficient to have two linearly independent lines in the sub-matrix \eqref{sub_matrix_obsv} to ensure the local weak observability of the system. 

\subsection{WRSM observability conditions}
The first two lines of \eqref{sub_matrix_obsv}, which come from the first derivatives of $i_{\alpha}$ and $i_{\beta}$, are studied. This choice is motivated by the fact that these currents are available for measurement in all synchronous machines, the rotor current (from which the third line of the matrix \eqref{sub_matrix_obsv} is calculated) does not exist in the case of PMSM and SyRM. Another reason comes from the physics of the machine: $i_f$ is a DC signal, whereas both $i_{\alpha}$ and $i_{\beta}$ are AC signals, so it is more convenient for physical interpretation to take AC signals together.

Symbolic math software is used to evaluate the determinant $\Delta_y$ of the sub-matrix composed of the first two lines of \eqref{sub_matrix_obsv}. In order to make the interpretation of this determinant easier, $\alpha \beta$ currents are expressed as functions of $dq$ currents\footnote{$dq$ currents are the machine currents in the rotating reference frame, which rotates at the rotor electrical speed (see Fig. \ref{wrsm_ab}). The machine equations in this reference frame are derived using the following Park transformation given by \eqref{Park_a} and \eqref{Park_b}.} using the Park transformation:
\begin{eqnarray}
i_\alpha &=& i_d \cos \theta - i_q \sin \theta \label{Park_a} \label{park1}\\
i_\beta &=& i_d \sin \theta + i_q \cos \theta \label{Park_b} \label{park2}
\end{eqnarray}

Finally, the determinant has the following form:
\begin{eqnarray}
\Delta_y &=& \mathcal{D}\omega + \mathcal{N}
\label{delta_wrsm}
\end{eqnarray}
where
\begin{eqnarray}
\mathcal{D} &=&  \frac{1}{L_D L_q} \left[
\left(L_\delta i_d + M_f i_f \right)^2 + L_\Delta L_\delta i_q^2
\right] \label{D_wrsm}\\
\mathcal{N} &=& \frac{L_\Delta}{L_D L_q} \left[
\left(L_\delta \frac{di_d}{dt} + M_f \frac{di_f}{dt} \right) i_q \right. \label{N_wrsm} \\
&& ~~~~~~~~~~~~~~~ - \left. \left(L_\delta i_d + M_f i_f \right) \frac{di_q}{dt}
\right] \nonumber
\end{eqnarray}
with
\begin{eqnarray}
L_\delta = L_d - L_q~;~
L_\Delta = L_\delta - \frac{M_f^2 }{L_f}~;~
L_D = L_d - \frac{M_f^2 }{L_f}
\end{eqnarray}
\begin{figure}
\centering
\includegraphics[scale=0.45]{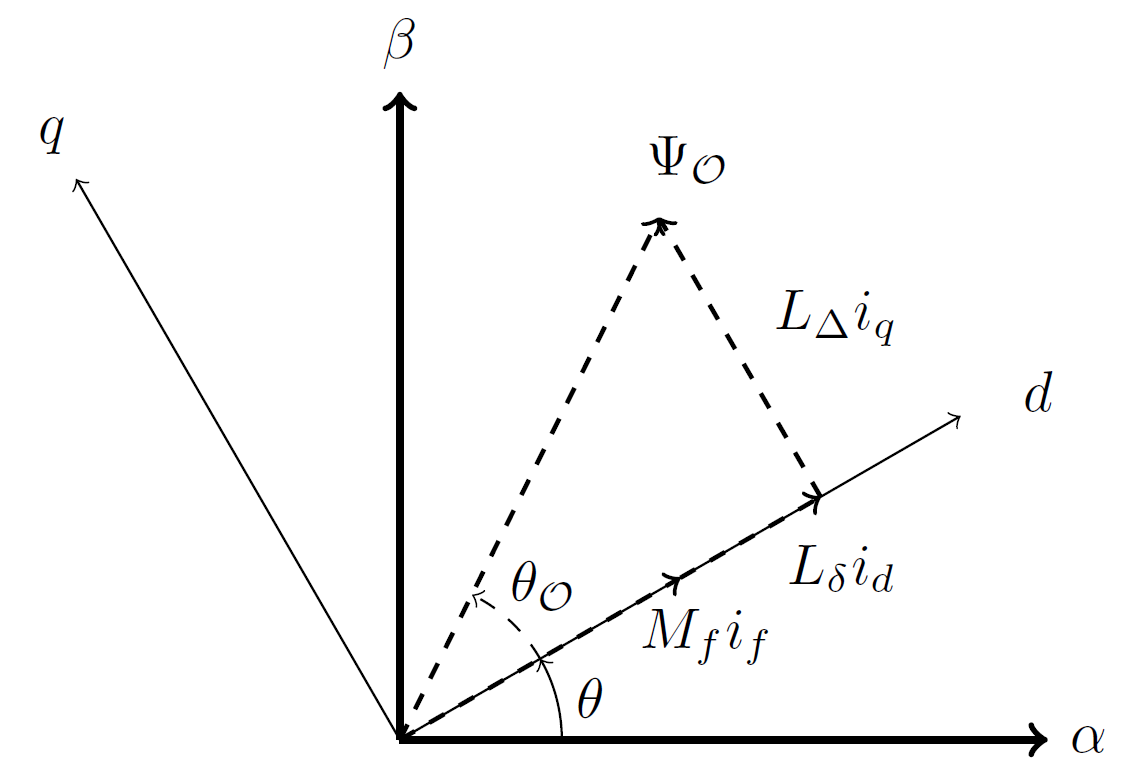}
\caption{Vector diagram of the WRSM showing stator reference frame (thick), rotor reference frame, and the observability vector (dashed)}
\label{wrsm_vector}
\end{figure}

The observability condition $\Delta_y \neq 0$ implies:

\begin{eqnarray}
\omega \neq \frac{\left(L_\delta i_d + M_f i_f \right) L_\Delta \frac{di_q}{dt} - 
\left(L_\delta \frac{di_d}{dt} + M_f \frac{di_f}{dt} \right) L_\Delta i_q}{\left(L_\delta i_d + M_f i_f \right)^2 + L_\Delta L_\delta i_q^2}
\end{eqnarray}
The above equation can be written as:
\begin{eqnarray}
\omega &\neq& \frac{\left(L_\delta i_d + M_f i_f \right)^2 + L_\Delta^2 i_q^2}{\left(L_\delta i_d + M_f i_f \right)^2 + L_\Delta L_\delta i_q^2} \times\\
&&\frac{\left(L_\delta i_d + M_f i_f \right) L_\Delta \frac{di_q}{dt} - 
\left(L_\delta \frac{di_d}{dt} + M_f \frac{di_f}{dt} \right) L_\Delta i_q}{\left(L_\delta i_d + M_f i_f \right)^2 + L_\Delta^2 i_q^2} \nonumber
\end{eqnarray}
then
\begin{eqnarray}
\omega &\neq& \frac{\left(L_\delta i_d + M_f i_f \right)^2 + L_\Delta^2 i_q^2}{\left(L_\delta i_d + M_f i_f \right)^2 + L_\Delta L_\delta i_q^2} \times \nonumber\\
&& ~~~~~~~~~~~~~~\frac{d}{dt} \arctan \left( \frac{L_\Delta i_q}{L_\delta i_d + M_f i_f} \right) 
\end{eqnarray}
The following approximation can be adopted\footnote{This approximation does not affect the observability conditions at standstill where $\omega = 0$ and currents are nonzero.}:
\begin{eqnarray}
\frac{\left(L_\delta i_d + M_f i_f \right)^2 + L_\Delta^2 i_q^2}{\left(L_\delta i_d + M_f i_f \right)^2 + L_\Delta L_\delta i_q^2} \approx 1
\label{assum}
\end{eqnarray}

Thus, the WRSM observability condition can be formulated as:
\begin{eqnarray}
\omega \neq \frac{d}{dt} \arctan \left( \frac{L_\Delta i_q}{L_\delta i_d + M_f i_f} \right) 
\label{obsv_vect_0}
\end{eqnarray}
It can be seen that the above equation describes a vector, which will be called the \textit{observability vector} and denoted $\Psi_\mathcal{O}$ (Fig. \ref{wrsm_vector}), that has the following components in the $dq$ reference frame:
\begin{eqnarray}
\Psi_{\mathcal{O}d} &=& L_\delta i_d + M_f i_f\\
\Psi_{\mathcal{O}q} &=& L_\Delta i_q
\end{eqnarray}
The condition \eqref{obsv_vect_0} becomes:
\begin{eqnarray}
\omega \neq \frac{d}{dt} \theta_{\mathcal{O}}
\label{obsv_vect_1}
\end{eqnarray}
where $\theta_{\mathcal{O}}$ is the phase of the vector $\Psi_\mathcal{O}$ in the rotor ($dq$) reference frame (see Fig. \ref{wrsm_vector}).

Finally, the following WRSM observability condition can be stated: the local observability of a WRSM is guaranteed if the rotational velocity of the observability vector with respect to the rotor is different from the electrical velocity of the rotor with respect to the stator. Therefore, at standstill, the observability vector should rotate and not be fixed. It turns out that the $d-$axis component of the observability vector is nothing but the \textit{active flux}, introduced by Boldea \textit{et al.} in \cite{boldea08af}, which is, by definition, the torque producing flux.

Obviously, if the (nonzero) currents $i_d$, $i_q$, and $i_f$ are constant at standstill, then the currents $i_\alpha$ and $i_\beta$ are also constant (this is straightforward from the equations \eqref{park1} and \eqref{park2}). In this case, the determinant \eqref{delta_wrsm} is equal to zero, and the observability condition is not fulfilled. To overcome this situation, we propose to inject a high-frequency (HF) current in the rotor winding in a way to make $i_f$ variable, so that the observability vector ``vibrates'' at standstill, and the observability condition \eqref{obsv_vect_1} is fulfilled. In practice, this technique can be useful for the starting of the machine, then, during the machine operation, the HF current will be injected only when the rotor estimated speed is near zero, in order to ensure the observability.

\subsection{Other SMs observability conditions}
The other SMs observability conditions can be derived from the previous results, taking into consideration the adequate equations of section \ref{subsec:ss_sms} for each machine. In addition, the following substitutions should be made for the PMSM and SyRM:
\begin{equation}
L_D = L_d ~~~;~~~ L_\Delta = L_\delta
\end{equation}
which means that the approximation \eqref{assum} is an equality for these machines.

The observability conditions interpretation can be generalized using the \textit{observability vector} concept. For instance, the observability vector of the SPMSM is equivalent to the rotor permanent magnet flux vector, then the only case where the observability is not guaranteed is the standstill (for further remarks on PMSM observability refer to \cite{koteich_sled_15}). Furthermore, the observability vector of the SyRM is aligned with the stator current space vector.

\section{Illustrative Simulations}
\label{simu}
\begin{figure}[!t]
\begin{center}
\includegraphics[scale=1.2]{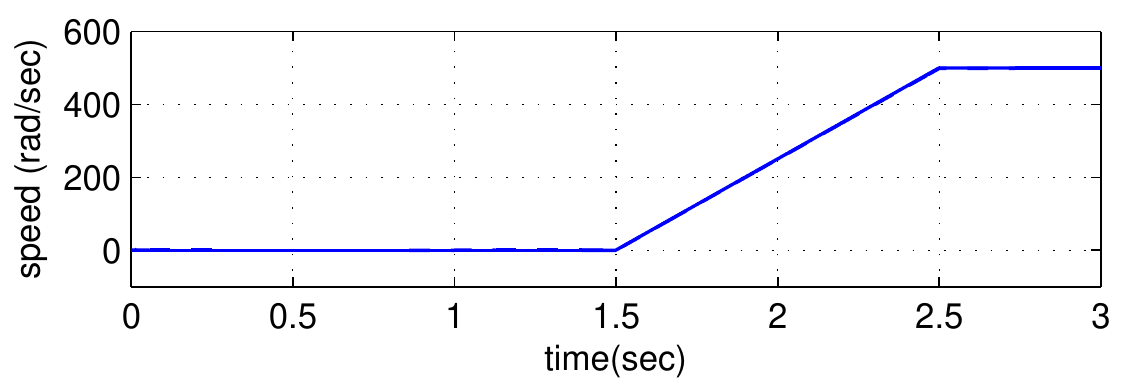}
\end{center}
\caption{Rotor speed profile}
\label{speed_profile}
\end{figure}

The present section is aimed at illustrating the previous observability analysis using numerical simulation. For this purpose, an extended Kalman filter (EKF) is designed. In order to make the study of some critical situations easier, the following operation mode is installed: the rotor position is considered to be driven by an external mechanical system, which imposes the following speed profile (Fig. \ref{speed_profile}): zero speed during $1.5~sec$, then a constant angular acceleration of $500~rd/s^{2}$ during one second, then the speed is fixed at $500~rd/s$. Stator and rotor currents are regulated, using standard proportional-integral (PI) controllers, to fit with the following set-points:
\begin{eqnarray}
\begin{array}{c c c c c}
i_d^* = 4~A & ; & i_q^* = 15~A & ; & i_f^* = 4~A
\end{array}
\end{eqnarray}
Table \ref{param} shows the machine parameters.

\subsection{Extended Kalman Filter}
The EKF algorithm is described below:
\subsubsection{Model linearization}
\begin{eqnarray}
A_k = \left.\frac{\partial f(x,u)}{\partial x}\right|_{x_k, u_k};~~~
C_k = \left.\frac{\partial h(x)}{\partial x}\right|_{x_k} 
\end{eqnarray}
\subsubsection{Prediction}
\begin{eqnarray}
\hat{x}_{k+1/k} &=& \hat{x}_{k/k} + T_s f(\hat{x}_{k/k},u_k)\\
P_{k+1/k} &=& P_k + T_s(A_k P_k + P_k A_k^T) + Q_k
\end{eqnarray}
\subsubsection{Gain}
\begin{equation}
K_k = P_{k+1/k} C_k^T(C_k P_{k+1/k} C_k^T + R_k)^{-1}
\end{equation}
\subsubsection{Innovation}
\begin{eqnarray}
\hat{x}_{k+1/k+1} &=& \hat{x}_{k+1/k} + K_k(y - h(\hat{x}_{k+1/k}))\\
P_{k+1/k+1} &=& P_{k+1/k} - K_k C_k P_{k+1/k}
\end{eqnarray}
where $T_s$ is the sampling period.
\subsubsection{Tuning} EKF tuning is done by the choice of covariance matrices $Q_k$ and $R_k$, using \emph{trial and error} method:
\begin{eqnarray}
Q_k = \left[
\begin{matrix}
\mathbb{I}_{3 \times 3} & \mathbb{O}_{3 \times 1} & \mathbb{O}_{3 \times 1}\\
\mathbb{O}_{1 \times 3} & 200 & 0\\
\mathbb{O}_{1 \times 3} & 0 & 5
\end{matrix}
\right]~;~
R_k = \mathbb{I}_{3 \times 3}
\end{eqnarray}

\subsection{HF current injection}
\begin{figure}[!t]
\begin{center}
\includegraphics[scale=1.2]{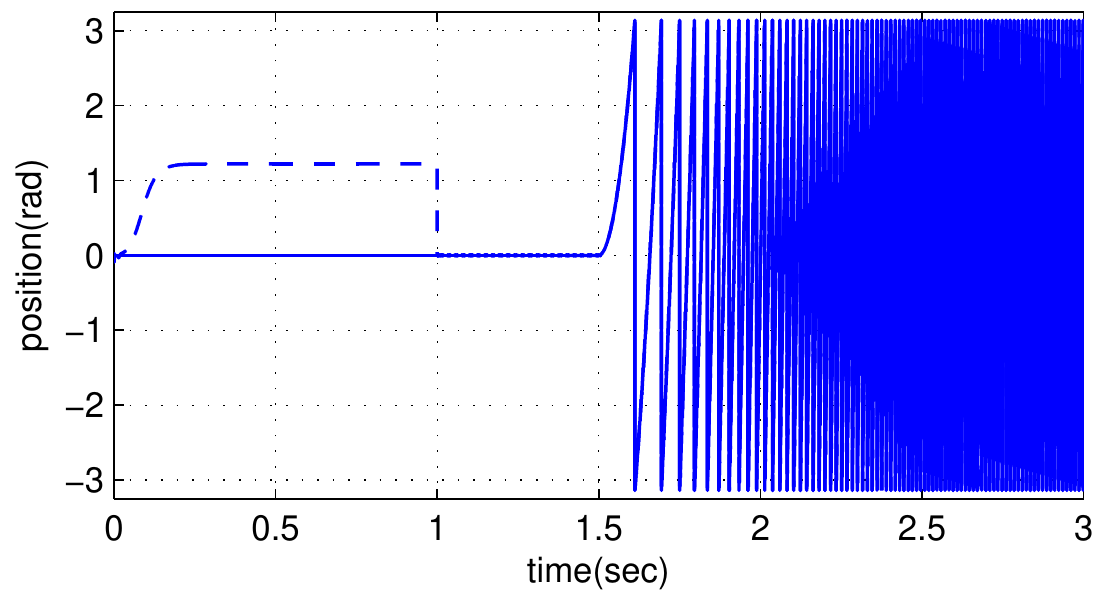}
\end{center}
\caption{Estimated rotor position (dashed) compared to the real position}
\label{position}
\end{figure}

\begin{figure}[!t]
\begin{center}
\includegraphics[scale=1.2]{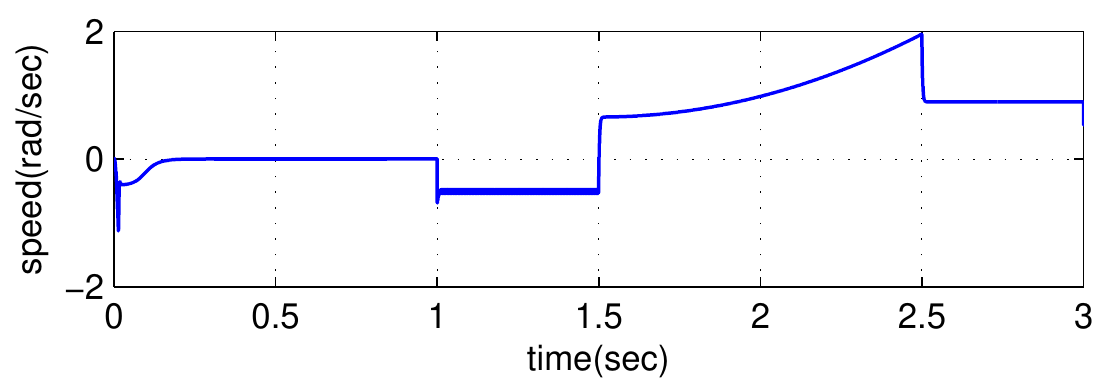}
\end{center}
\caption{Rotor speed estimation error}
\label{speed}
\end{figure}

The following HF current is added to the rotor current $i_f$ during the time interval $[1~s., 1.5~s.]$:
\begin{eqnarray}
i_{f_{HF}} = I_{f_{HF}} \sin \omega_{HF}t = 0.5 \sin 2 \pi 10^3 t ~~ A
\end{eqnarray}
Fig.~\ref{position} shows the real and estimated rotor angular positions. It is obvious that, at standstill, the EKF does not converge to the correct value of $\theta$ until the HF current is injected. For nonzero speeds, there is no position estimation problem. The speed estimation error is shown in Fig. \ref{speed}; the error slightly increases with the HF injection, but it remains reasonable.


The choice of the injected signal amplitude $I_{f_{HF}}$ and its angular frequency $\omega_{HF}$ has to be done taking into consideration some practical aspects: very high frequencies generate more losses in the magnetic circuit, however, low frequencies generate undesired vibration in the motor. High amplitude HF current generates both more losses and more vibration, whereas low amplitude (and very high frequencies) might be filtered by the rotor electrical inertia without any effect on the observability.

\section{Conclusions}
The concept of \textit{observability vector} is introduced in this paper. The observability analysis of sensorless synchronous machine drives shows that the local observability can be guaranteed if the rotational speed of the observability vector with respect to the rotor is different from the electrical angular speed of the rotor with respect to the stator. 

Based on the above analysis, a high-frequency current injection technique is proposed for the wound-rotor synchronous machine sensorless control; it consists of injecting an HF alternating current in the rotor windings when the rotor speed is near zero, which makes the observability vector vibrate around its position.

The unified approach adopted in this paper can be useful in finding similar solutions for the other synchronous dives.

\begin{table}[!t]
\renewcommand{\arraystretch}{1.3}
\label{table}
\caption{Parameters of the WRSM Used in simulation}
\centering
\begin{tabular}{|l|c|}
\hline
Parameters & Value [Unit]\\
\hline
\hline
Number of pole pairs (${p}$) & 2\\
Stator resistance $R_s$ & $0.01$ [$\Omega$]\\
Rotor resistance $R_f$ & $6.5$ [$\Omega$]\\
Direct inductance $L_d$ & $0.8$ [$mH$]\\
Quadratic inductance $L_q$ & $0.7$ [$mH$]\\
Rotor inductance $L_f$ & $0.85$ [$H$]\\
\hline
\end{tabular}
\label{param}
\end{table}
\bibliographystyle{ieeetr}
\bibliography{biblio_koteich}

\begin{thebibliography}{10}

\bibitem{emadi14}
B.~Bilgin and A.~Emadi, ``Electric motors in electrified transportation: A step
  toward achieving a sustainable and highly efficient transportation system,''
  {\em Power Electronics Magazine, IEEE}, vol.~1, pp.~10--17, June 2014.

\bibitem{mea12}
X.~Roboam, B.~Sareni, and A.~Andrade, ``More electricity in the air: Toward
  optimized electrical networks embedded in more-electrical aircraft,'' {\em
  Industrial Electronics Magazine, IEEE}, vol.~6, pp.~6--17, Dec 2012.

\bibitem{chiasson}
J.~Chiasson, {\em Modeling and high performance control of electric machines},
  vol.~26.
\newblock John Wiley \& Sons, 2005.

\bibitem{boldea08}
I.~Boldea, ``Control issues in adjustable speed drives,'' {\em Industrial
  Electronics Magazine, IEEE}, vol.~2, pp.~32--50, Sept 2008.

\bibitem{vas}
P.~Vas, {\em Sensorless vector and direct torque control}.
\newblock Monographs in electrical and electronic engineering, Oxford: Oxford
  University Press, 1998.

\bibitem{holtz06}
J.~Holtz, ``Sensorless control of induction machines - with or without signal
  injection?,'' {\em Industrial Electronics, IEEE Transactions on}, vol.~53,
  pp.~7--30, Feb 2005.

\bibitem{koteich13}
M.~Koteich, T.~Le~Moing, A.~Janot, and F.~Defay, ``A real-time observer for
  uav's brushless motors,'' in {\em Electronics, Control, Measurement, Signals
  and their application to Mechatronics (ECMSM), 2013 IEEE 11th International
  Workshop of}, pp.~1--5, June 2013.

\bibitem{betin14}
F.~Betin, G.-A. Capolino, D.~Casadei, B.~Kawkabani, R.~Bojoi, L.~Harnefors,
  E.~Levi, L.~Parsa, and B.~Fahimi, ``Trends in electrical machines control:
  Samples for classical, sensorless, and fault-tolerant techniques,'' {\em
  Industrial Electronics Magazine, IEEE}, vol.~8, pp.~43--55, June 2014.

\bibitem{hermann77}
R.~Hermann and A.~J. Krener, ``Nonlinear controllability and observability,''
  {\em Automatic Control, IEEE Transactions on}, vol.~22, pp.~728--740, Oct
  1977.

\bibitem{besancon}
G.~Besancon, {\em Nonlinear observers and applications}.
\newblock Lecture Notes in Control and Information Sciences, Verlag/Heidelberg,
  New-York/Berlin: Springer, 2007.

\bibitem{canudas00}
C.~de~Wit, A.~Youssef, J.~P. Barbot, P.~Martin, and F.~Malrait, ``Observability
  conditions of induction motors at low frequencies,'' in {\em Decision and
  Control, 2000. Proceedings of the 39th IEEE Conference on}, vol.~3,
  pp.~2044--2049 vol.3, 2000.

\bibitem{rojas}
S.~Ibarra-Rojas, J.~Moreno, and G.~Espinosa-Pérez, ``Global observability
  analysis of sensorless induction motors,'' {\em Automatica}, vol.~40, no.~6,
  pp.~1079 -- 1085, 2004.

\bibitem{chiasson06}
M.~Li, J.~Chiasson, M.~Bodson, and L.~Tolbert, ``A differential-algebraic
  approach to speed estimation in an induction motor,'' {\em Automatic Control,
  IEEE Transactions on}, vol.~51, pp.~1172--1177, July 2006.

\bibitem{zaltni10}
D.~Zaltni, M.~Ghanes, J.~P. Barbot, and M.-N. Abdelkrim, ``Synchronous motor
  observability study and an improved zero-speed position estimation design,''
  in {\em Decision and Control (CDC), 2010 49th IEEE Conference on},
  pp.~5074--5079, Dec 2010.

\bibitem{ezzat10}
M.~Ezzat, J.~de~Leon, N.~Gonzalez, and A.~Glumineau, ``Observer-controller
  scheme using high order sliding mode techniques for sensorless speed control
  of permanent magnet synchronous motor,'' in {\em Decision and Control (CDC),
  2010 49th IEEE Conference on}, pp.~4012--4017, Dec 2010.

\bibitem{abry11}
F.~Abry, A.~Zgorski, X.~Lin-Shi, and J.-M. Retif, ``Sensorless position control
  for spmsm at zero speed and acceleration,'' in {\em Power Electronics and
  Applications (EPE 2011), Proceedings of the 2011-14th European Conference
  on}, pp.~1--9, Aug 2011.

\bibitem{vaclavek13}
P.~Vaclavek, P.~Blaha, and I.~Herman, ``Ac drive observability analysis,'' {\em
  Industrial Electronics, IEEE Transactions on}, vol.~60, pp.~3047--3059, Aug
  2013.

\bibitem{park29}
R.~Park, ``Two-reaction theory of synchronous machines generalized method of
  analysis-part i,'' {\em American Institute of Electrical Engineers,
  Transactions of the}, vol.~48, pp.~716--727, July 1929.

\bibitem{park33}
R.~Park, ``Two-reaction theory of synchronous machines-ii,'' {\em American
  Institute of Electrical Engineers, Transactions of the}, vol.~52,
  pp.~352--354, June 1933.

\bibitem{boldea08af}
I.~Boldea, M.~Paicu, and G.~Andreescu, ``Active flux concept for
  motion-sensorless unified ac drives,'' {\em Power Electronics, IEEE
  Transactions on}, vol.~23, pp.~2612--2618, Sept 2008.

\bibitem{koteich_sled_15}
M.~Koteich, A.~Maloum, G.~Duc, and G.~Sandou, ``Permanent magnet synchronous
  drives observability analysis for motion-sensorless control,'' in {\em
  Sensorless Control for Electrical Drives (SLED), 2015 IEEE Symposium on},
  June 2015.

\end{thebibliography}
\end{document}